\documentclass[11pt]{article}
\usepackage{latexsym}
\usepackage{amsmath, amsthm}
\usepackage[all,dvips,ps]{xy}
\usepackage{amsfonts}
\usepackage{mathrsfs}

\newtheorem{lem}{Lemma}
\newtheorem{thm}{Theorem}
\newtheorem{cor}{Corollary}
\newtheorem{exa}{Example}
\newtheorem{ass}{Assumption}
\newtheorem{defn}{Definition}

\newcommand{\Rs}{\mathbb{R}}

\newcommand{\Sn}{{\cal S}^n }
\newcommand{\Snp}{{\cal S}^n_+ }

\newcommand{\snn}{{\cal S}^{n-1} }
\newcommand{\snnp}{{\cal S}^{n-1}_+ }
\newcommand{\Dn}{{\cal D}^{n} }

\newcommand{\mm}{\overline{m} }

\newcommand{\tl}{\langle}
\newcommand{\tr}{\rangle}
\newcommand{\tra}{\mbox{trace} }
\newcommand{\diag}{\mbox{diag} }
\newcommand{\rank}{\mbox{rank} }

\newcommand{\N}{\mathscr{N} }

\newcommand{\U}{\mathcal{U} }
\newcommand{\W}{\mathcal{W} }

\newcommand{\E}{\mathcal{E} }

\newcommand{\Fs}{\mathscr{F} }
\newcommand{\Ls}{\mathscr{L} }
\newcommand{\Xs}{\mathscr{X} }
\newcommand{\Ys}{\mathscr{Y} }
\newcommand{\T}{\mathscr{T}}
\newcommand{\K}{\mathscr{K} }
\newcommand{\TV}{\mathscr{T}_V }
\newcommand{\KV}{\mathscr{K}_V }

\newcommand{\bz}{{\bf 0 }}

\newcommand{\bpr}{{\bf Proof.} \hspace{1 em}}
\newcommand{\epr}{ \\ \hspace*{4.5in} $\Box$ }
\newcommand{\beq}{ \begin{equation} }
\newcommand{\eeq}{ \end{equation} }
\newcommand{\bt}{ \begin{tabular} }
\newcommand{\et}{ \end{tabular} }
\begin{document}

\bibliographystyle{plain}
\title{Universal Rigidity of Bar Frameworks via the Geometry of Spectrahedra
 \thanks{Research supported by the Natural Sciences and Engineering
         Research Council of Canada.} }
\vspace{0.3in}
        \author{ A. Y. Alfakih
  \thanks{E-mail: alfakih@uwindsor.ca}
  \\
          Department of Mathematics and Statistics \\
          University of Windsor \\
          Windsor, Ontario N9B 3P4 \\
          Canada }

\date{\today}
\maketitle


\vspace{0.1in}

\begin{abstract}
A bar framework $(G,p)$ in dimension $r$ is a graph $G$ whose nodes are points $p^1,\ldots,p^n$ in
$\Rs^r$ and whose edges are line segments between pairs of these points.
Two frameworks $(G,p)$ and $(G,q)$ are equivalent if
each edge of $(G,p)$ has the same (Euclidean) length as the corresponding edge of $(G,q)$.
A pair of non-adjacent vertices $i$ and $j$ of $(G,p)$ is universally linked if $||p^i-p^j||$ = $||q^i-q^j||$
in every framework $(G,q)$ that is equivalent to $(G,p)$. Framework $(G,p)$ is universally rigid {\em iff}
every pair of non-adjacent vertices of $(G,p)$ is universally linked.

In this paper, we present a unified treatment of the universal rigidity problem based on the geometry of spectrahedra.
A spectrahedron is the intersection of the positive semidefinite cone with an affine space.
This treatment makes it possible to tie together some known, yet scattered, results and to derive new ones.
Among the new results presented in this paper are: (i)
A sufficient condition for a given pair of
non-adjacent vertices of $(G,p)$ to be universally linked. (ii)
A new, weaker, sufficient condition for a framework $(G,p)$ to be universally rigid thus strengthening the existing
known condition. An interpretation of this new condition in terms of the Strong Arnold Property and transversal
intersection is also presented.
\end{abstract}

\section{Introduction and Main Results}
\label{int}

A {\em bar framework} (or a framework for short) $(G,p)$ in dimension $r$
is a simple graph $G$ whose vertices are
points $p^1,\ldots,p^n$ in $\Rs^r$ and whose edges are line segments between pairs of
these points. $(G,p)$ is $r$-dimensional if the points $p^1,\ldots,p^n$ affinely span $\Rs^r$.
To avoid trivialities we assume that $G$ is connected and not complete.
Each framework $(G,p)$ generates a Euclidean distance matrix $D_p = (d_{ij}) = ||p^i - p^j||^2$,
where $||.||$ denotes the Euclidean norm. Let $H$ denote the adjacency matrix of $G$ and let
$A \circ B$ denote the Hadamard, or element-wise, product of matrices $A$ and $B$, i.e., $(A\circ B)_{ij} = A_{ij} B_{ij}$.
Two frameworks $(G,p)$ and $(G,q)$ are
{\em equivalent} if each edge of $(G,p)$ has the same Euclidean length as the corresponding edge of $(G,q)$, i.e., if
\[ H \circ D_p = H \circ D_q.
\]

A pair of non-adjacent vertices $i$ and $j$ of $(G,p)$ is {\em universally linked}
if $||p^i-p^j||$=$||q^i-q^j||$ in every framework $(G,q)$ that is equivalent to $(G,p)$.
Framework $(G,p)$ is {\em universally rigid} if and only if every non-adjacent pair of vertices
of $(G,p)$ is universally linked; i.e.,
$(G,p)$ is universally rigid {\em iff} $H \circ D_p = H \circ D_q$ implies that $D_p=D_q$.

The notions of dimensional rigidity and affine motions turned out to be very useful in the study of universal rigidity.
An $r$-dimensional framework $(G,p)$ is {\em dimensionally rigid} if no
$s$-dimensional framework $(G,q)$, for any $s \geq r+1$, is equivalent to $(G,p)$.
We say that a framework $(G,p)$ has an {\em affine motion} if there exists a framework $(G,q)$ such that:
(i) $(G,q)$ is equivalent to $(G,p)$, (ii) $D_p \neq D_q$,  and (iii) $q^i=Ap^i+ b$ for all $i=1,\ldots,n$,
for some matrix $A$ and vector $b$.

The following theorem gives a characterization of universal rigidity in terms of dimensional rigidity and affine motions.
\begin{thm}[Alfakih \cite{alf07a}] \label{thmuda}
A bar framework $(G,p)$ is universally rigid if and only if it is dimensionally rigid and has no affine motion.
\end{thm}

Stress matrices play a key role in the problems of universal and dimensional rigidities.
An {\em equilibrium stress} (or simply a stress) of $(G,p)$ is a real-valued function
$\omega$ on $E(G)$, the edge set of $G$, such that
\beq \label{defw}
\sum_{j:\{i,j\} \in E(G)} \omega_{ij} (p^i - p^j) = \bz \mbox{ for all } i=1,\ldots,n.
\eeq
Let $\overline{E}(G)$ denote the set of missing edges of $G$, i.e.,
\[
\overline{E}(G)= \{ \{i,j\}: i \neq j , \{i,j\} \not \in E(G) \},
\]
and let $\omega =(\omega_{ij})$ be a stress of $(G,p)$. Then the
$n \times n$ symmetric matrix $\Omega$ where
\beq \label{defO}
\Omega_{ij} = \left\{ \begin{array}{ll} -\omega_{ij} & \mbox{if } \{i,j\} \in E(G), \\
                        0   & \mbox{if }  \{i,j\}  \in \overline{E}(G), \\
                   {\displaystyle \sum_{k:\{i,k\} \in E(G)} \omega_{ik}} & \mbox{if } i=j,
                   \end{array} \right.
\eeq
is called the {\em stress matrix} associated with $\omega$, or a stress matrix
of $(G,p)$. It is not hard to show that the maximum rank of any stress matrix of
an $r$-dimensional framework on $n$ vertices is $n-r-1$.

The following theorem provides a sufficient condition for the
universal rigidity of a given framework.

\begin{thm}[Connelly \cite{con82}, Alfakih \cite{alf07a}]
\label{thmsuff1}
Let $(G,p)$ be an $r$-dimensional framework on $n$ vertices in $\Rs^r$, for some $r
\leq n-2$. If the
following two conditions holds:
\begin{enumerate}
\item $(G,p)$ admits a positive semidefinite stress matrix $\Omega$ of rank
$n-r-1$,
\item $(G,p)$ has no affine motion,
\end{enumerate}
then $(G,p)$ is universally rigid.
\end{thm}

Sufficient conditions for a framework to have no affine motion were presented in \cite{con05,ay13,an13}.
Proofs that the reverse of Theorem \ref{thmsuff1} holds were given in \cite{gt14} if $(G,p)$ is generic,
in \cite{alf11} if $G$ is a chordal graph, and in
\cite{aty13} if $G$ is an $(r+1)$-lateration graph.

In this paper, we present a unified treatment of the universal rigidity problem based on the geometry of spectrahedra.
A spectrahedron is the intersection of the cone of positive semidefinite matrices with an affine space \cite{rg95}.
In addition to obtaining new results,
this treatment makes it possible to tie together some known, yet scattered, results.
In particular, the geometric approach developed in this paper leads to new simple proofs of Theorems \ref{thmuda}
and \ref{thmsuff1}.

The following sufficient condition for a given pair of
non-adjacent vertices of $(G,p)$ to be universally linked is the first new main result of this paper.
Throughout the paper, let  $E^{ij}$ denote the $n \times n$ symmetric matrix with 1's in the $(i,j)$th and $(j,i)$th entries
and 0's elsewhere.

\begin{thm}\label{thmmain1}
Let $(G,p)$ be an $r$-dimensional framework on $n$ vertices in $\Rs^r$, $r \leq (n-2)$, and let
$\{k,l\} \in \overline{E}(G)$. Let $\Omega$ be a nonzero positive semidefinite
stress matrix of $(G,p)$. If the following condition holds:
\beq \label{eqnmain1}
 \mbox{There does not exist } y_{kl}  \neq 0 : \Omega (\sum_{\{i,j\} \in \overline{E}(G)} y_{ij} E^{ij}) = \bz,
\eeq
then the pair $\{k,l\}$ is universally linked.
\end{thm}
The proof of Theorem \ref{thmmain1} is given at the end of Section \ref{sub}.
Unfortunately, as Example \ref{ex2} below shows, Condition (\ref{eqnmain1}) is not necessary.
Moreover, it is worth pointing out that Condition (\ref{eqnmain1}) is equivalent to:
\beq \label{eqnmain11}
\mbox{There does not exist } y_{kl}  \neq 0 : \sum_{\{i,j\} \in \overline{E}(G)} y_{ij} E^{ij} \W = \bz,
\eeq
where $\W$ is the matrix whose columns form a basis of the column space of $\Omega$.

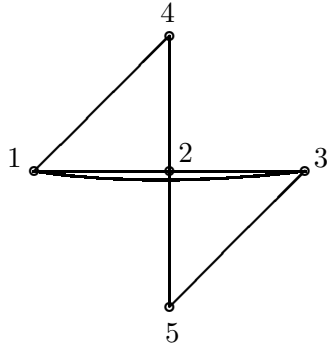
\begin{figure}
\thicklines
\setlength{\unitlength}{0.6mm}
\begin{picture}(80,80)(-70,-25)

\put(30,10){\circle{2}}
\put(60,10){\circle{2}}
\put(90,10){\circle{2}}
\put(60,40){\circle{2}}
\put(60,-20){\circle{2}}

\put(24,11){$1$}
\put(62,12){$2$}
\put(92,11){$3$}
\put(58,43){$4$}
\put(59,-28){$5$}

\put(30,10){\line(1,1){30}}
\put(30,10){\line(1,0){30}}
\put(60,10){\line(1,0){30}}
\put(60,10){\line(0,1){30}}
\put(60,10){\line(0,-1){30}}
\put(60,-20){\line(1,1){30}}
\qbezier(30,10)(56,6)(90,10)

\end{picture}
\caption{The $2$-dimensional frameworks in $\Rs^2$ of Example \ref{ex1}.
The edge $\{1,3\}$ is drawn as an arc to make edges $\{1,2\}$ and $\{2,3\}$
visible. The pairs of non-adjacent vertices $\{1,5\}$ and $\{3,4\}$ are universally linked.
Obviously, the pair of non-adjacent vertices $\{4,5\}$ is not universally linked. }
\label{f1}
\end{figure}

\begin{exa} \label{ex1}
Consider the $2$-dimensional framework $(G,p)$ on 5 vertices depicted in Figure \ref{f1}.
This framework admits a nonzero positive semidefinite stress matrix whose column space
is spanned by
$\W = \left[ \begin{array}{r} 1 \\ -2 \\ 1 \\ 0 \\ 0 \end{array} \right]$. In this case,
the system of equations:
 \[ \sum_{\{i,j\} \in \overline{E}(G)} y_{ij} E^{ij} \W =
 \left[ \begin{array}{ccccc} 0 & 0 & 0 & 0 & y_{15} \\ 0 & 0 & 0 & 0 & 0 \\
                               0 & 0 & 0 & y_{34} & 0 \\ 0 & 0 & y_{34} & 0 & y_{45} \\
                                y_{15} & 0 & 0 & y_{45} & 0 \end{array} \right]
 \left[ \begin{array}{r} 1 \\ -2 \\ 1 \\ 0 \\ 0 \end{array} \right] =
 \left[ \begin{array}{r} 0 \\ 0 \\ 0 \\ 0 \\ 0 \end{array} \right]
 \]
 has the solution
 $y_{15}=y_{34}=0$ while $y_{45}$ is free.
 Therefore the non-adjacent pairs $\{1,5\}$ and $\{3,4\}$ are universally linked. Note that the pair
 $\{4,5\}$ is not universally linked since the framework can fold along the edge $\{1,3\}$.
 \end{exa}

The second new main result of this paper is the following theorem which presents
a weaker sufficient condition for universal rigidity than that in Theorem \ref{thmsuff1}.
In fact, this new condition does not assume that the rank of stress matrix $\Omega$ is
equal to $n-r-1$.

\begin{thm} \label{thmmain2}
Let $(G,p)$ be an $r$-dimensional framework on $n$ vertices in $\Rs^r$, $r \leq (n-2)$.
Let $\Omega$ be a nonzero positive semidefinite
stress matrix of $(G,p)$. If the following condition holds:
\beq \label{eqnmain2}
 \mbox{There does not exist } y=(y_{ij})  \neq \bz:  \Omega(\sum_{\{i,j\} \in \overline{E}(G)} y_{ij} E^{ij}) = \bz,
\eeq
then $(G,p)$ is universally rigid.
\end{thm}
Even though it follows as an immediate corollary of Theorem \ref{thmmain1}, we present
a direct proof of Theorem \ref{thmmain2} at the end of Section \ref{sub}. Unfortunately, as Example \ref{ex2} below shows,
Condition~(\ref{eqnmain2}) is not necessary. It is worth noting that if stress matrix $\Omega$ in Theorem \ref{thmmain2} has rank = $n-r-1$,
then, as Corollary \ref{corend} below shows, Condition~(\ref{eqnmain2}) is in fact equivalent to the assertion that
framework $(G,p)$ has no affine motions.
Thus, Theorem \ref{thmsuff1} is a special case of Theorem \ref{thmmain2}.
Finally, we point out that
Condition~(\ref{eqnmain2}) is equivalent to:
\beq \label{eqnmain22}
\mbox{There does not exist } y=(y_{ij})  \neq \bz : \sum_{\{i,j\} \in \overline{E}(G)} y_{ij} E^{ij} \W = \bz,
\eeq
where $\W$ is the matrix whose columns form a basis of the column space of $\Omega$.

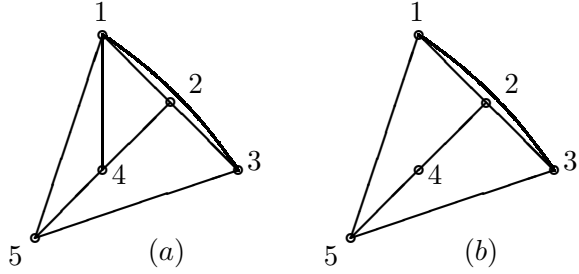
\begin{figure}
\thicklines
\setlength{\unitlength}{0.6mm}
\begin{picture}(70,70)(-70,-25)

\put(0,30){\circle{2}}
\put(15,15){\circle{2}}
\put(30,0){\circle{2}}
\put(0,0){\circle{2}}
\put(-15,-15){\circle{2}}

\put(-2,33){$1$}
\put(19,17){$2$}
\put(32,0){$3$}
\put(2,-3){$4$}
\put(-21,-21){$5$}

\put(-15,-15){\line(1,3){15}}
\put(-15,-15){\line(1,1){15}}
\put(-15,-15){\line(3,1){45}}
\put(0,0){\line(1,1){15}}
\put(0,30){\line(1,-1){30}}
\put(0,0){\line(0,1){30}}
\qbezier(0,30)(18,18)(30,0)

\put(10,-20){$(a)$}

\put(70,30){\circle{2}}
\put(85,15){\circle{2}}
\put(100,0){\circle{2}}
\put(70,0){\circle{2}}
\put(55,-15){\circle{2}}

\put(68,33){$1$}
\put(89,17){$2$}
\put(102,0){$3$}
\put(72,-3){$4$}
\put(49,-21){$5$}

\put(55,-15){\line(1,3){15}}
\put(55,-15){\line(1,1){15}}
\put(55,-15){\line(3,1){45}}
\put(70,0){\line(1,1){15}}
\put(70,30){\line(1,-1){30}}
\qbezier(70,30)(88,18)(100,0)

\put(80,-20){$(b)$}

\end{picture}
\caption{The two universally rigid $2$-dimensional frameworks in $\Rs^2$ of Example \ref{ex2}.
The edge $\{1,3\}$ is drawn as an arc to make edges $\{1,2\}$ and $\{2,3\}$
visible. Both frameworks have a positive semidefinite stress matrix $\Omega$ of rank 1,  }
\label{f2}
\end{figure}

\begin{exa} \label{ex2}
Consider the $2$-dimensional framework $(G,p)$ on 5 vertices depicted in Figure \ref{f2} (a).
This framework has the same stress matrix $\Omega$ of Example \ref{ex1}.
Thus even though this framework is obviously universally rigid, it admits
a positive semidefinite stress matrix of rank $1 < 2=n-r-1$. Thus Theorem \ref{thmsuff1} does not apply in this case.
On the other hand, the solution of the system of equations $(y_{25} E^{25} + y_{34} E^{34})\W=\bz$ is
$y_{25}=y_{34}=0$. Thus  Condition (\ref{eqnmain2}) holds in this case.

To show that Condition (\ref{eqnmain2})
is not necessary, consider the framework in Figure \ref{f2} (b). In this case,
$(y_{14} E^{14}+ y_{25} E^{25} + y_{34} E^{34})\W=\bz$ has the solution:
$y_{25}=0$ and $y_{14}=- y_{34}=1$.
Thus Condition (\ref{eqnmain2}) does not hold even though this framework is universally rigid.
\end{exa}

We conclude this section with the following interpretation of Condition (\ref{eqnmain2}) in terms
of the Strong Arnold Property of matrices and the notion of transverse intersection of manifolds. Note that
these two notions, in our case, are equivalent to the notion of non-degeneracy in semidefinite programming \cite{aho97,lv14}.

Given a graph $G$ on $n$ vertices, Let $A$ be a matrix in $\Sn$,
the space of $n \times n$ symmetric matrices, such that $A_{ij}=0$ for all $\{i,j\} \in \overline{E}(G)$.
$A$ is said to satisfy the {\em Strong Arnold Property (SAP)} if $Y=\bz$ is the only matrix in $\Sn$ that satisfies:
(i) $Y_{ij} = 0$ if $i=j$ or if $\{i,j\} \in E(G)$, (ii) $AY=\bz$.
Therefore, Equation (\ref{eqnmain2}) is equivalent to the assertion that stress matrix $\Omega$ satisfies
the (SAP).

Another notion which is equivalent to the SAP, in our case, is that of transversal intersection.
Let $\Omega$ be a stress matrix of rank $k$, and let $\U$ and $\W$ be the matrices whose
columns form, respectively, orthonormal bases of the null space and the column space (or the range) of $\Omega$.
Let $\Sn_k=\{ A \in \Sn : \mbox{ rank } A = k\}$. Then the tangent space  of $\Sn_k$
at $\Omega$ is given by
\beq
{\cal T}_{\Omega}  =  \{ A \in \Sn : A = [ \W \; \U]
\left[ \begin{array}{cc} \Phi_1 & \Phi_2 \\
                        \Phi_2^T & \bz \end{array} \right]
\left[ \begin{array}{c} \W^T \\ {\U}^T \end{array} \right] \},
\eeq
 where $\Phi_1$ is a symmetric matrix of order $k$ and $\Phi_2$ is
$k \times (n-k)$.
Further, let
\[ {\cal L} = \{  A \in \Sn: A_{ij} = \tra(A E^{ij}) = 0 \mbox{ whenever } \{i,j\} \in \overline{E}(G) \}. \]
Thus $\Omega \in \Sn_k \cap {\cal L}$.
We say that $\Sn_k$ intersects ${\cal L}$ {\em  transversally} at $\Omega$ if
\beq \label{eqnti1}
{\cal T}_{\Omega} + {\cal L} = \Sn.
\eeq
But Equation (\ref{eqnti1}) is equivalent to
\beq \label{eqnti2}
{\cal T}_{\Omega}^{\perp} \cap {\cal L}^{\perp} = \{\bz\},
\eeq
where
\[
{\cal T}_{\Omega}^{\perp} = \{ A \in \Sn : A = \U \Phi {\U}^T \}.
\]
and
\[ {\cal L}^{\perp} = \mbox{ span } (\{E^{ij}: \{i,j\} \in \overline{E}(G) \}). \]
In other words, Equation (\ref{eqnti2}) is equivalent to:
\[
\mbox{There does not exist } y=(y_{ij}) \neq \bz:  \U \Phi \U^T = \sum_{\{i,j\} \in \overline{E}(G)} y_{ij} E^{ij}.
\]
Observe that $y=\bz$ iff $\Phi=\bz$ since $E^{ij}$'s are linearly independent.
Therefore, Equation (\ref{eqnmain2}) is equivalent to the assertion that the set of symmetric matrices
of rank $k$ intersects the subspace ${\cal L}$ transversally at $\Omega$.

\subsection{Notation}

For easy reference, the notation used throughout the paper are collected here.
$\Sn$ denotes the space of $n \times n$ symmetric matrices equipped with the usual
inner product $\tl A, B \tr$ = trace$\,(AB)$.
The positive semi-definiteness (definiteness) of a real symmetric matrix $A$
is denoted by $A \succeq 0$ ($A \succ 0$) and the cone of $n \times n$ symmetric
positive semidefinite matrices is denoted by $\Snp$.
We denote by $E^{ij}$ the $n \times n$ symmetric matrix with 1's in the $(i,j)$th and $(j,i)$th entries and 0's elsewhere,
and by $e$ the vector of all 1's in $\Rs^n$.
We denote the identity matrix of order $n$ by $I_n$.
$\bz$ denotes the zero matrix or zero vector of appropriate dimension.
$E(G)$ denotes the set of edges of a simple graph $G$ and
$\overline{E}(G)$ denotes the set of the missing edges $G$, i.e.,
$\overline{E}(G)=\{ \{i,j\} : i \neq j, \{i,j\} \not \in E(G) \}$.
The cardinality of $\overline{E}(G)$, i.e., the number of missing edges of $G$, is denoted by
$\mm$. For a matrix $A$, $\diag(A)$ denotes the vector consisting of the diagonal entries of $A$.
Finally, $\N(A)$ denote the null space of matrix $A$.

\section{Preliminaries}
\label{pre}

This section develops the necessary mathematical background needed in this paper.
In particular, we review basic definitions and results concerning Euclidean distance
matrices, projected Gram matrices and Gale matrices.

\subsection{Euclidean Distance Matrices}
An $n \times n$ symmetric matrix $D=(d_{ij})$ is called a {\em Euclidean distance matrix (EDM)} if there exist points
$p^1, \ldots, p^n$ in some Euclidean space such that $ d_{ij}=||p^i-p^j||^2$ for all $i,j=1,\ldots,n$.
The dimension of the affine span
of $p^1,\ldots,p^n$ is called the {\em embedding dimension} of $D$.
Let $e$ denote the vector of all 1's in $\Rs^n$ and let $e^{\perp}$ denote
the orthogonal complement of $e$ in $\Rs^n$, i.e.,
\[
e^{\perp} = \{ x \in \Rs^n: e^T x = 0 \}.
\]
Let $V$ be the $n \times (n-1)$ matrix whose columns form an orthonormal basis of $e^{\perp}$, i.e.,
$V$ satisfies:
\beq \label{defV}
V^Te=\bz \mbox{ and } V^TV=I_{n-1}.
\eeq
Define \cite{cri88} the linear operator
\beq \label{defT}
\T (D) := - \frac{1}{2} J D J,
\eeq
where
\[ J= VV^T = I - \frac{1}{n} ee^T \]
is the orthogonal projection onto $e^{\perp}$.
Then it is well known \cite{cri88,gow85,sch35,yh38} that a symmetric matrix $D$ whose diagonal entries are all
0's is EDM with edmedding dimension $r$ if and only if
\beq \label{eqnsch}
\T(D) \succeq \bz \mbox{ and } \rank(\T(D)) = r.
\eeq
Let $\Dn$ denote the set of EDMs of order $n$, then it immediately follows from (\ref{eqnsch}) that
$\Dn$ is a convex cone. Moreover,
\beq
\T(\Dn) = \{ B \in \Snp : Be = 0 \},
\eeq
where $\Snp$ denotes the cone of $n \times n$ symmetric positive semidefinite matrices.
Define the two subspaces of $\Sn$,  the space of $n \times n$ symmetric matrices.
\begin{eqnarray*}
S_H & = & \{ A \in \Sn : \diag(A) = \bz \} , \\
S_C & = & \{ B \in \Sn : Be = \bz \},
\end{eqnarray*}
where $\diag(A)$ denotes the vector consisting of the diagonal entries of $A$.
Let \cite{cri88}
\beq \label{defK}
\K (B) := \diag(B)e^T + e (\diag(B))^T - 2 B.
\eeq
\begin{thm}[Gower \cite{gow85}]
$\K(S_C) = S_H$, $\T(S_H) = S_C$ and $\K_{|S_c}$ and $\T_{|S_H}$ are inverses of each other.
\end{thm}
Thus
\beq
\K( \{ B \in \Snp : Be = 0 \}) = \Dn.
\eeq
As will be shown below, the set $ \{ B \in \Snp: Be = 0 \}$ is a maximal proper face of $\Snp$.
Moreover, every face of $\Snp$ is isomorphic to a positive semidefinite cone of lower dimension. In particular,
\[  \{ B \in \Snp:  Be = 0 \} = \{ B= V X V^T : X \in \snnp \}.
\]
Therefore, define \cite{akw99} the two linear operators
\begin{eqnarray*}
\KV (X) & := & \K(VXV^T) = \diag(VXV^T)e^T + e (\diag(VXV^T))^T - 2 VXV^T, \\
\TV(D) & := & V^T \T (D)V = - \frac{1}{2} V^T D V.
\end{eqnarray*}
Hence,
\beq
\KV(\snnp) = \Dn \mbox{  and  } \TV(\Dn) = \snnp .
\eeq
Therefore, EDMs are in a one-to-one correspondence with positive semidefinite matrices of order $n-1$.

\subsection{The Configuration Spectrahedron of a Framework}
\label{conspec}

The Cayley configuration space of  a framework $(G,p)$ is the set of all possible distance values of the
missing edges of $(G,p)$ when $(G,p)$ is viewed as a linkage consisting of rigid bars (edges) and universal joints (vertices).
In this section we present a characterization of such configuration space.
This characterization plays a pivotal role in our approach. Our presentation here follows closely \cite{alf00}.

The {\em configuration matrix} $P$ of  an $r$-dimensional framework
$(G,p)$ on $n$ vertices in $\Rs^r$ is the $n\times r$ matrix whose $i$th row is $(p^i)^T$, i.e.,
\[
P=\left[ \begin{array}{c} (p^1)^T \\ \vdots \\ (p^n)^T \end{array} \right].
\]
Consequently, the Gram matrix of $(G,p)$ is  $PP^T$. We make the following assumption throughout this paper.
\begin{ass} \label{ass1}
In any framework $(G,p)$, the centroid of the points $p^1,\ldots,p^n$
coincides with the origin, i.e., $P^Te= \bz$,
\end{ass}

Let $D_p=(d_{ij})$ be the EDM generated by $(G,p)$, then
\begin{eqnarray*}
d_{ij} & = & ||p^i - p^j ||^2, \\
       & = & (p^i)^T p^i + (p^j)^T p^j - 2 (p^i)^T p^j, \\
       & = & (PP^T)_{ii} + (PP^T)_{jj} - 2 (PP^T)_{ij}.
\end{eqnarray*}
Therefore, $D_p = \K(PP^T)$ = $\KV(X)$, where
 \beq
 PP^T=VXV^T, \mbox{ and } X = V^TPP^TV.
 \eeq
$X$ is called the {\em projected Gram matrix} of $(G,p)$.
Observe that $X$ is a symmetric positive semidefinite matrix of order $n-1$ and of rank $r$.
Let $X'$ be the projected Gram matrix of framework $(G,p')$. Then $(G,p')$ is equivalent to $(G,p)$ if
and only if
\beq \label{eqnHKV}
H \circ (D_p - D_{p'}) = H \circ \KV(X - X') = \bz,
\eeq
where $H$ is the adjacency matrix of $G$. Recall that
$E^{ij}$ denotes the $n \times n$ symmetric matrix with 1's in the $(i,j)$th and $(j,i)$th entries and 0's elsewhere.
The set $\{E^{ij}: 1 \leq i < j \leq n\}$ is obviously linearly independent.  Let
\beq \label{defM}
M^{ij} = \TV(E^{ij}) = - \frac{1}{2} V^T E^{ij} V.
\eeq
\begin{lem} \label{lemMli}
The set $\{ M^{ij} : \{i,j\} \in \overline{E}(G) \}$ is a basis of the null space of $H \circ \KV$.
\end{lem}
\bpr
First we show that the set $\{ M^{ij} : 1 \leq i<j \leq n \}$ is linearly independent.
Assume that $\sum_{ i < j } \alpha_{ij} M^{ij}$ = $\bz$, then
$\KV(\TV (\sum_{i < j} \alpha_{ij} E^{ij})) = \bz$.
But $ E^{ij} \in S_H $ if $ i < j$. Thus, $\sum_{i < j} \alpha_{ij} E^{ij} = \bz$, which implies that
$\alpha_{ij}=0$ for all $ i < j$. Hence, the set
 $\{ M^{ij} : \{i,j\} \in \overline{E}(G) \}$ is linearly independent. Now
for all $\{i,j\} \in \overline{E}(G)$ we have $\KV(M^{ij}) = \KV(\TV(E^{ij})) = E^{ij}$.
Hence, $H \circ \KV(M^{ij})= H \circ E^{ij} = \bz$.
\epr

\begin{defn} \label{defF}
Let $X$ denote the projected Gram matrix of framework $(G,p)$.
Let $M^{ij}$ be the matrices defined in (\ref{defM}) and let
$\mm$ denote the cardinality of $\overline{E}(G)$.
Further, for $y=(y_{ij}) \in \Rs^{\mm}$, define
\beq
\Xs(y):=  X + \sum_{\{i,j\} \in \overline{E}(G)} y_{ij} M^{ij}.
\eeq
Then the spectrahedron
\beq
\Fs = \{ y   \in \Rs^{\mm} : \Xs(y) \succeq \bz \}
\eeq
is called the Cayley configuration spectrahedron of $(G,p)$.
\end{defn}
Clearly, $\Fs$ is a closed convex set in $\Rs^{\mm}$.
The following theorem justifies Definition~\ref{defF}.
\begin{thm} \label{thmjustF}
Let $\Fs$ be the Cayley configuration spectrahedron of an $r$-dimensional framework $(G,p)$ and
let $X'$  be the projected Gram matrix of an $s$-dimensional framework $(G,p')$.
Then $(G,p')$ is equivalent to $(G,p)$ if and only if
\[
X'= \Xs (y) \mbox{ for some } y \in \Fs \mbox{ with } \rank \; X' = s,
\]
in which case, $||{p'}^i - {p'}^j ||^2 = ||p^i-p^j||^2 + y_{ij}$
for all $\{i,j\} \in \overline{E}(G)$.
\end{thm}
\bpr
The first part of the theorem follows from Lemma \ref{lemMli} and Equation (\ref{eqnHKV}).
Now $||{p'}^i - {p'}^j ||^2$ = $(\KV(X'))_{ij}$ = $(\KV(X))_{ij}$ + $\sum_{\{k,l\} \in \overline{E}(G)} y_{kl} (\KV(M^{kl}))_{ij}$.
But $\KV(M^{kl})=E^{kl}$.
Therefore,
$ ||{p'}^i - {p'}^j ||^2$ = $||p^i-p^j||^2+ y_{ij}$ if $ \{i,j\} \in \overline{E}(G)$. Obviously,
 $ ||{p'}^i - {p'}^j ||^2$  = $||p^i-p^j||^2$ if $\{i,j\} \in E(G)$.
\epr

The following two theorems are immediate consequences of Theorem~\ref{thmjustF} and Definition~\ref{defF}.
\begin{thm} \label{thmul}
Let $\Fs$ be the Cayley configuration spectrahedron of framework $(G,p)$ and
let $\{k,l\} \in \overline{E}(G)$.
Then $\{k,l\}$ is universally linked if and only if $\Fs$ is contained in the subspace
$\{ y \in \Rs^{\mm}: y_{kl}=0 \}$ of $\Rs^{\mm}$.
\end{thm}
\begin{thm} \label{thmurF0}
Let $\Fs$ be the Cayley configuration spectrahedron of framework $(G,p)$.
Then $(G,p)$ is universally rigid if and only if $\Fs=\{\bz\}$.
\end{thm}

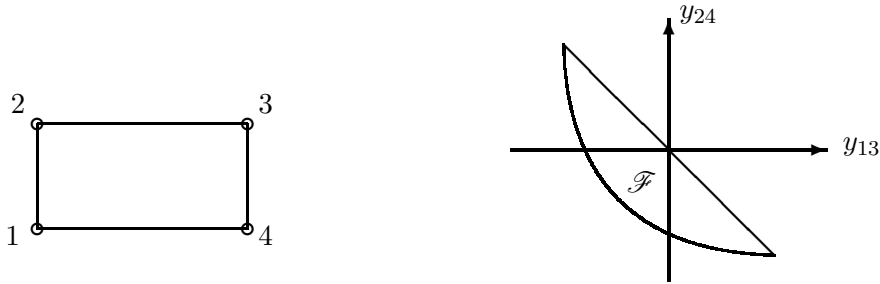
\begin{figure}
\thicklines
\setlength{\unitlength}{0.7mm}
\begin{picture}(80,80)(-20,-25)

\put(0,10){\circle{2}}
\put(40,10){\circle{2}}
\put(0,30){\circle{2}}
\put(40,30){\circle{2}}

\put(-6,7){$1$}
\put(42,7){$4$}
\put(-5,32){$2$}
\put(42,32){$3$}

\put(0,10){\line(1,0){40}}
\put(0,10){\line(0,1){20}}
\put(40,30){\line(-1,0){40}}
\put(40,30){\line(0,-1){20}}

\put(90,25){\vector(1,0){60}}
\put(120,0){\vector(0,1){50}}

\put(153,25){$y_{13}$}
\put(122,50){$y_{24}$}

\put(100,45){\line(1,-1){40}}
\put(112,17){$\Fs$}
\qbezier(100,45)(101,6)(140,5)

\end{picture}
\caption{The $2$-dimensional framework of Example \ref{ex3}
and its Cayley configuration space $\Fs$. }
\label{f3}
\end{figure}

We conclude this subsection with the following example of spectrahedron $\Fs$.
\begin{exa} \label{ex3}
Consider the $2$-dimensional framework $(G,p)$ in the plane depicted in Figure \ref{f3}, where
\[ p^1= \left[ \begin{array}{c} -1  \\ -0.5  \end{array} \right],
p^2= \left[ \begin{array}{c} -1  \\ 0.5  \end{array} \right],
p^3= \left[ \begin{array}{c} 1  \\ 0.5  \end{array} \right]
\mbox{ and } p^4= \left[ \begin{array}{c} 1  \\ -0.5  \end{array} \right].
\]
Then after straightforward but tedious calculation, we
get that the Cayley configuration space $\Fs$ of $(G,p)$ is given by
\[
\Fs=\{ y= \left[ \begin{array}{c} y_{13} \\ y_{24} \end{array} \right]
             : \frac{1}{2} \left[ \begin{array}{ccc} 2  & 2+y_{13} & -y_{24} \\
              2+y_{13}  & 10+2y_{13} & 8+y_{13} \\ -y_{24} & 8+y_{13} & 8  \end{array} \right] \succeq \bz \}.
\]
i.e.,
\[
\Fs=\{ y= \left[ \begin{array}{c} y_{13} \\ y_{24} \end{array} \right]
                             : y_{24} \leq - y_{13} \mbox{ and } y_{24} \geq - \frac{5 y_{13} + 16}{y_{13}+5} \}.
\]
Spectrahedron $\Fs$ is also depicted in Figure~\ref{f3}.
It is easy to see that $-4 \leq y_{13}, y_{24} \leq 4$.
Therefore, in any framework $(G,p')$ equivalent to $(G,p)$, $1 \leq ||{p'}^1 - {p'}^3||^2 \leq 9$
and $1 \leq ||{p'}^2 - {p'}^4||^2 \leq 9$
 since $||p^1 - p^3||^2=||p^2 - p^4||^2=5$.
\end{exa}

\subsection{Stress Matrices and Gale Matrices}

Stress matrices, which play a key role in the problem of universal rigidity, are closely related
to Gale matrices (or Gale transform \cite{gal56,gru67}). Gale matrices arise naturally when dealing with
the null space of projected Gram matrices. In this subsection, we review the basic results on stress matrices and
Gale matrices.

 A Gale matrix of an $r$-dimensional framework $(G,p)$ on $n$ vertices in $\Rs^r$, $r \leq n-2$, is any
$n \times (n-r-1)$ matrix $Z$ whose columns form a basis of
$\N( \left[ \begin{array}{c} P^T \\ e^T \end{array} \right])$, where $\N$ denotes the null space and
$P$ is the configuration matrix
of $(G,p)$. Observe that the matrix $[P \;\; e]$ has full column rank since the points
$p^1,\ldots,p^n$ affinely span $\Rs^r$.

It easily follows from the definition of a stress matrix in (\ref{defO}) that
a symmetric matrix $\Omega$ is a stress matrix of $(G,p)$ if and only if
\[
\Omega P =\bz, \Omega e =\bz \mbox{ and } \Omega_{ij} = 0 \mbox{ for all } \{i,j\} \in \overline{E}(G).
\]

\begin{lem}[\cite{alf10}] \label{lemOZ}
Let $Z$ be a Gale matrix of an $r$-dimensional framework $(G,p)$ on $n$ vertices.
Then $\Omega$ is a stress matrix of $(G,p)$ if and only if
\[
\Omega = Z \Psi Z^T \mbox{ with $(Z\Psi Z^T)_{ij} = 0$ for all } \{i,j\} \in \overline{E}(G),
\]
where $\Psi$ is a symmetric matrix of order $ n-r-1$.
\end{lem}
An immediate consequence of Lemma \ref{lemOZ} is that rank $\Omega$ = rank $\Psi \leq n-r -1$.
The relationship between Gale matrices and the null space of projected Gram matrices is established in
the following lemma.

\begin{lem}[\cite{alf07a}] \label{lemVU}
Let $X$  be the projected Gram matrix of framework $(G,p)$ and
let $U$ be the matrix whose columns form an orthonormal basis of $\N(X)$. Then
\[
VU \mbox{ is a Gale matrix of } (G,p),
\]
where $V$ is as defined in (\ref{defV}).
\end{lem}

\section{The Geometry of the Configuration Spectrahedron $\Fs$}

We begin this section by investigating the facial structure of $\Fs$ where we present a characterization of dimensional rigidity
in terms of the relative interior of $\Fs$. Then we focus on
the affine hulls of the faces of $\Fs$ and their connection with affine motions of the framework.

\subsection{Facial Structure of $\Fs$}

A set $S$ in $\Rs^n$ is affine if the line joining any two points of
$S$ lies entirely in $S$, i.e., for any two points $x^1$ and $x^2$ in $S$, the point
$x = \lambda x^1 + (1-\lambda) x^2$ is in $S$ for all scalar $\lambda$.
Every affine set in $\Rs^n$ is parallel to a linear subspace of $\Rs^n$.
The dimension of an affine set is equal to the dimension of the linear subspace parallel to it.
The {\em affine hull} of a set $S$ in $\Rs^n$, denoted by aff($S$), is the smallest affine set containing $S$.
The dimension of a set $S$ in $\Rs^n$ is equal to the dimension of aff($S$).
A set consisting of a single point is its own affine hull as well as its own relative interior, i.e.,
aff($\{\hat{x}\})$ =  $\{\hat{x}\}$ = relint($\{\hat{x}\})$.

A spectrahedron is the set formed by the intersection of the positive semidefinite cone with an affine set.
Let $F$ be a nonempty convex subset of spectrahedron ${\cal F}$. Then
$F$ is said to be a {\em face} of ${\cal F}$ if for every
$y$ in $F$ such that $y=\lambda x+ (1-\lambda) z$ for some $x$ and $z$ in ${\cal F}$ and some $\lambda: 0 < \lambda < 1$, it
follows that $x$ and $z$ are both in  $F$. Note that $\emptyset$ and ${\cal F}$, itself, are faces of ${\cal F}$. These two faces
${\cal F}$ are called {\em improper faces} while the other faces are called {\em proper faces}.
Let $y \in {\cal F}$ and let face($y$) denote the unique minimal face of ${\cal F}$ containing $y$ in its relative interior.
The following characterization of minimal faces, stated in terms of $\Fs$, is well known \cite{bc75,rg95,pat00}.
\begin{thm}
Let $y \in \Fs$, the Cayley configuration spectrahedron of framework $(G,p)$. Then
\begin{eqnarray*}
\mbox{face}(y) & = & \{ x \in \Fs: \N(\Xs(y)) \subseteq \N(\Xs(x)) \}, \\
\mbox{relint(face}(y)) & = & \{ x \in \Fs: \N(\Xs(y)) = \N(\Xs(x))\}.
\end{eqnarray*}
\end{thm}
Note that the faces of $\Fs$ are characterized by their relative interior.
Let $S$ be a convex subset of $\Fs$. The minimal face of $\Fs$ containing $S$, denoted by face($S$),
is the face of $\Fs$ which intersects the relative interior of $S$.

\begin{thm} \label{thmrankrelint}
Let $S$ be a convex subset of $\Fs$, the Cayley configuration spectrahedron of framework $(G,p)$.
Let $\hat{y}$ be a point in $S$, then the following statements are equivalent:
\begin{enumerate}
\item rank $\Xs(\hat{y}) \geq $ rank $\Xs(y)$ for all $y \in S$,
\item face($\hat{y}$) =  face($S$),
\item  $\hat{y} \in$ relint($S$).
\end{enumerate}
\end{thm}
\bpr
(1 $\implies$ 2)

Since $\hat{y} \in S$, it follows that face($\hat{y}) \subseteq $ face($S$). Now let $y$ be any point
in $S$ and let $x = \lambda \hat{y} + (1-\lambda) y$ for some $\lambda: 0 < \lambda <1$.
Then $x \in S$ since $S$ is convex. Moreover,
$\Xs(x) = \lambda \Xs(\hat{y}) + (1-\lambda) \Xs(y)$. Therefore, $\N(\Xs(x)) \subseteq \N(\Xs(\hat{y}))$ and
$\N(\Xs(x)) \subseteq \N(\Xs(y))$, and hence, rank $\Xs(x) \geq$ rank  $\Xs(\hat{y})$. But by assumption
rank $\Xs(x) \leq$ rank  $\Xs(\hat{y})$. Therefore, rank $\Xs(x)$ = rank  $\Xs(\hat{y})$ and consequently
$\N(\Xs(x)) = \N(\Xs(\hat{y}))$. Hence, $\N(\Xs(\hat{y})) \subseteq  \N(\Xs(y))$. Therefore,
$y \in $ face($\hat{y}$) and thus $S \subset $ face($\hat{y}$). Hence, face($S) \subseteq $ face($\hat{y}$).

(2 $\implies$ 3)

Let face($\hat{y}$) = face ($S$) and assume that $\hat{y} \not \in $ relint($S$). Then
$\hat{y}$ is in the relative boundary of $S$. Hence,
there exists a hyperplane ${\cal H}$ containing $\hat{y}$ but not $S$ \cite{roc70}. Thus
${\cal H} \cap $ face($S$) is a face of $\Fs$, containing $\hat{y}$, smaller than
face($S$) = face($\hat{y}$) which contradicts the definition of face($\hat{y}$).

(3 $\implies$ 1)

Assume that $\hat{y} \in$ relint($S$) and let $y \in S$. Then there exists $x$ in $S$ and $\lambda: 0 < \lambda < 1$
such that $\hat{y}=\lambda y + (1-\lambda) x$. Hence, $\N(\Xs(\hat{y})) \subseteq  \N(\Xs(y))$
since $\Xs(\hat{y}) = \lambda \Xs(y) + (1-\lambda) \Xs(x)$. Therefore,
rank $\Xs(\hat{y}) \geq$ rank  $\Xs(y)$.
\epr

The following characterization of dimensional rigidity of $(G,p)$ is an immediate consequence
of Theorem \ref{thmrankrelint}.
\begin{thm} \label{thmdimrig0}
Let $\Fs$ be the Cayley configuration spectrahedron of framework $(G,p)$.
Then $(G,p)$ is dimensionally rigid if and only if $ \bz \in $ relint($\Fs$).
\end{thm}
\bpr
Set $S=\Fs$ and $\hat{y}= \bz$ in Theorem \ref{thmrankrelint}.
\epr

\subsection{Affine motion of $(G,p)$}

We focus in this subsection on the affine motions of a framework $(G,p)$ and their connection to
the affine hulls of the faces of $\Fs$.

\begin{lem} \label{lemaffF}
Let $x \in \Fs$, the Cayley configuration spectrahedron of framework $(G,p)$.
Let $U$ be the matrix whose columns  form
an orthonormal basis for $\N(\Xs(x))$. Then the affine hull of face($x$) is given by
\beq \label{eqnam}
\mbox{aff(face}(x))=  \{ y \in \Rs^{\mm}: \Xs(y)U = \bz \}.
\eeq
\end{lem}
\bpr
Let $\Ls=\{ y \in \Rs^{\mm}: \Xs(y)U  = \bz \}$ and let $z$ be a point in face($x$). Further,
let $y=\lambda x + (1-\lambda) z $ for some scalar $\lambda$. Then
$\Xs(y)U= \lambda \Xs(x) U + (1-\lambda) \Xs(z) U = \bz$. Hence, $y \in \Ls$ and thus
aff(face($x$)) $ \subseteq \Ls$.

On the other hand, if $\Ls = \{x\}$, then we are done. Therefore,  let $y \in \Ls$, $y \neq x$ and
let $W$ be the matrix whose columns  form
an orthonormal basis for the column space of $\Xs(x)$.
Since $x$ is in relint(face($x$)), it follows that $W^T \Xs(x)W \succ \bz$. Hence,
there exists $t>0$ such that $W^T \Xs(x)W - t ( W^T\Xs(y)W - W^T \Xs(x)W) \succeq \bz$.
Let $z=x-t(y-x)$. Then $W^T \Xs(z) W \succeq \bz$ and  $\Xs(z)U=\bz$. Let $Q=[U \; W]$, then
$Q^T \Xs(z) Q \succeq \bz$ and hence, $\Xs(z) \succeq \bz$. Therefore, $z \in $ face($x$). Moreover,
since $y=(1+1/t)x- z/t$, it follows that $y$ belongs to the affine hull of face($x$) and hence
$\Ls \subseteq $ aff(face($x$)).
\epr

\begin{lem} \label{lemaffm}
Let $\Fs$ be the Cayley configuration spectrahedron of framework $(G,p)$ and let
$Z$ be a Gale matrix of $(G,p)$. Then the following statements are equivalent:
\begin{enumerate}
\item $(G,p)$ does not have an affine motion,
\item $y = \bz$ is the only solution of the equation $ \sum_{\{i,j\} \in \overline{E}(G)}y_{ij} V^T E^{ij} Z = \bz$,
\item aff(face($\bz$))= $\{\bz\}$.
\end{enumerate}
\end{lem}
\bpr The equivalence between statements 1 and 2 was first proved in \cite{alf10}. Now let $X$ be the projected
Gram matrix of $(G,p)$ and let $U$ be the matrix whose columns form and orthonormal basis of $\N(X)$.
Then it follows from Lemma \ref{lemaffF} and Lemma \ref{lemVU} that
\[ \mbox{aff(face($\bz$))} = \{ y \in \Rs^{\mm}: \sum_{\{i,j\} \in \overline{E}(G)}y_{ij} V^T E^{ij} Z = \bz\}.
\]
\epr

Theorem \ref{thmdimrig0} and Lemma \ref{lemaffm} lead to the following immediate proof of
Theorem~\ref{thmuda}. Assume that framework $(G,p)$ is dimensionally rigid and has no affine motion. Then
$\bz \in $ relint($\Fs$) and aff(face($\bz$)) = $\{\bz\}$. But in this case face($\bz$)= face($\Fs$) = $\Fs$. Therefore,
aff($\Fs$) = $\{\bz\}$ and hence, $\Fs=\{\bz\}$. Therefore, by Theorem \ref{thmurF0}, framework $(G,p)$ is universally rigid.
Finally, we remark that the framework of Example \ref{ex3} has an affine motion since 
aff(face($\bz$)) = $\{ y \in \Rs^2 : y_{13} + y_{24}=0 \}$ as can be seen in Figure \ref{f3}.

\section{The Subspace Containing $\Fs$}
\label{sub}

Our investigation in the previous section did make use of stress matrices. In fact, as well be shown below, if a framework
$(G,p)$ admits a nonzero positive semidefinite stress matrix $\Omega$, then $\Fs$, the Cayley configuration space
of $(G,p)$, is contained in a subspace of $\Rs^m$ determined by $\Omega$. An investigation of this subspace is given
in this section. First we begin with
following well-known Farkas lemma on the positive semidefinite cone which will be used to characterize the existence
of a nonzero positive semidefinite stress matrix.

\begin{lem} \label{farkas}
Let $A^0,A^1,\ldots,A^m$ be given matrices in $\snn$.
Then exactly one of the following two statements hold:
\begin{enumerate}
\item There exists $x=(x_{ij}) \in \Rs^m$ such that $A^0+ x_1 A^1+ \cdots + x_m A^m \succ \bz$,
\item There exists $Y \succeq \bz$, $Y \neq \bz$ such that $ \tra(A^0 Y) \leq 0$ and $\tra(Y A^i) = 0$ for all $i=1,\ldots,m$.
\end{enumerate}
\end{lem}

\begin{thm} \label{farkas2}
Let $\Fs$ be the Cayley configuration spectrahedron of framework $(G,p)$.
Then exactly one of the following two statements hold:
\begin{enumerate}
\item There exists $y \in \Fs$ such that $\Xs(y) \succ \bz$.
\item There exists $\Ys \succeq \bz, \neq \bz$ such that $X\Ys=\bz$ and $\tra(\Ys M^{ij}) = 0$ for all $\{i,j\} \in \overline{E}(G)$.
\end{enumerate}
\end{thm}
\bpr
Since $X \succeq \bz$ and $\Ys \succeq \bz$, $\tra(X \Ys) \leq 0$ implies that $X \Ys=\bz$. Thus, the result follows from
Lemma \ref{farkas}.
\epr

Theorem \ref{farkas2} provides an immediate proof to the following result.

\begin{thm}[Alfakih \cite{alf11}] \label{thmfarapp}
Let $(G,p)$ be a framework. Then exactly one of the following two statements hold:
\begin{enumerate}
\item $\exists$ $(n-1)$-dimensional framework $(G,p')$ that is equivalent to $(G,p)$.
\item $(G,p)$ admits a stress matrix $\Omega \succeq \bz, \neq \bz$.
\end{enumerate}
\end{thm}
\bpr Let $X$ be the projected Gram matrix of $(G,p)$ and let $U$ be the matrix whose
columns form an orthonormal basis of $\N(X)$. Then $X\Ys=\bz$ iff $\Ys=U \Psi U^T$ for some
symmetric matrix $\Psi$.
Let $\Omega=V \Ys V^T$ and let $\{i,j\} \in \overline{E}(G)$.
Then it follows from Lemma \ref{lemVU} that
$\Omega = VU \Psi U^T V^T= Z \Psi Z^T$, where $Z$ is a Gale matrix of $(G,p)$.
Moreover, $\Omega e= \bz$ and $\Omega_{ij}$ = $\tra(\Omega E^{ij})$ = $-2 \; \tra(\Ys M^{ij})$ = 0. Hence,
it follows from Lemma \ref{lemOZ} that $\Omega$ is a stress matrix of $(G,p)$. Moreover,
$\Ys \succeq \bz$ iff $\Omega \succeq \bz$ since $\Omega e= \bz$.
Therefore, Statement 2 of Lemma \ref{farkas2} is equivalent to the existence of a nonzero
positive semidefinite stress matrix.
\epr

The framework on $4$ vertices considered in Example \ref{ex3}, obviously, has an equivalent $3$-dimensional framework.
Equally obvious is the fact that
this framework has no nonzero positive semidefinite stress matrix.

\begin{lem} \label{lemOX}
Let $\Fs$ be the Cayley configuration space of framework $(G,p)$ and let $\Omega$ be
a nonzero positive semidefinite stress matrix of $(G,p)$.  Then
\[
 \Omega V \Xs(y) V^T = \bz \mbox{ for all } y \in \Fs.
\]
\end{lem}
\bpr
Let $y \in \Fs$, then
\[ \tra(\Omega V \Xs(y) V^T) = \tra(\Omega VXV^T) + \sum_{\{i,j\} \in \overline{E}(G)} y_{ij} \tra(\Omega E^{ij}) = 0.
\]
Therefore, $\Omega V \Xs(y) V^T = \bz$ since $\Omega \succeq \bz$ and $V \Xs(y) V^T \succeq \bz$.
\epr

Lemma \ref{lemOX} provides a simple proof of Theorem \ref{thmsuff1} since if
$\Omega$ is a positive semidefinite stress matrix of $(G,p)$
of rank $n-r-1$, then rank $V \Xs(y)V^T $ = rank $\Xs(y)  \leq r$ for all $y \in \Fs$.
Hence, rank $X$= rank $\Xs(\bz) \geq $ rank $\Xs(y)$ for all $y \in \Fs$,
and therefore, $(G,p)$ is dimensionally rigid.
Another immediate consequence of Lemma \ref{lemOX} is the following corollary.

\begin{cor}
Let $\Omega$ be a nonzero positive semidefinite stress matrix
of $(G,p)$.  If $(G,p')$ is a framework that is equivalent to $(G,p)$, then
$\Omega$ is a stress matrix of $(G,p')$.
\end{cor}

\bpr
Let $P'$ be the configuration matrix of $(G,p')$,  then $P'{P'}^T = V \Xs(y)V^T$
for some $y \in \Fs$. Thus it follows from Lemma \ref{lemOX}
that $\Omega P' {P'}^T$ = $\Omega V \Xs(y) V^T =\bz$, and hence $\Omega P'=\bz$.
\epr
\begin{lem} \label{lemxi}
Let $\Omega$ be a stress matrix of framework $(G,p)$. Then
\[
\{y : \Omega(\sum_{\{i,j\} \in \overline{E}(G)} y_{ij} E^{ij}V) = \bz\} \mbox{ is equal to }
\{y : \Omega(\sum_{\{i,j\} \in \overline{E}(G)} y_{ij}E^{ij}) = \bz\}
\]
\end{lem}
\bpr
Let $\E(y)= \sum_{\{i,j\} \in \overline{E}(G)} y_{ij} E^{ij}$ and let
$\Omega \E(y) V = \bz$. Then
$\Omega \E(y) = \xi e^T$ for some $\xi \in \Rs^n$. To complete the proof it suffices to show that $\xi=\bz$.
To this end, for all $i=1,\ldots,n$, we have
\begin{eqnarray*}
\xi_i & = &  \sum_{k=1}^n \Omega_{ik} \E(y)_{ki} \\
      & = & \Omega_{ii} \E(y)_{ii} + \sum_{k: \{i,k\} \in E(G) } \Omega_{ik} \E(y)_{ki} +
                                    \sum_{k: \{i,k\} \in \overline{E}(G) } \Omega_{ik} \E(y)_{ki} \\
      &= & 0,
\end{eqnarray*}
since $\E(y)_{ii}=0$, $\E(y)_{ki}=0$ for all $\{i,k\} \in E(G)$ and $\Omega_{ik}=0$ for all
$\{i,k\} \in \overline{E}(G)$.
\epr

The following theorem which characterizes the subspace containing $\Fs$ is the basis of the proofs of our main results.

\begin{thm}\label{thmproof}
Let $(G,p)$ be a framework and let $\Omega$ be a nonzero positive semidefinite stress matrix
of $(G,p)$. Then $\Fs$, the Cayley configuration spectrahedron of $(G,p)$,
is contained in the subspace
\beq \label{eqnYX0}
\{ y \in \Rs^{\mm}: \Omega ( \sum_{\{i,j\} \in \overline{E}(G)} y_{ij} E^{ij}) = \bz \}.
\eeq
\end{thm}
\bpr
Let $y \in \Fs$. Then it follows from Lemma \ref{lemOX} that $\Omega V \Xs(y) V^T=\bz$.
Thus,
\[
\Omega V \Xs(y) = \Omega ( \sum_{\{i,j\} \in \overline{E}(G)} y_{ij} E^{ij}V ) = \bz
\]
since $\Omega VV^T = \Omega (I - ee^T/n)$ = $\Omega$. The result follows from Lemma \ref{lemxi}.
\epr

\noindent{\bf Proof of Theorem \ref{thmmain1}}
Let $\Fs$ be the Cayley configuration space of framework $(G,p)$ and
suppose that Condition \ref{eqnmain1} of Theorem \ref{thmmain1} holds. Then it follows from
Theorem~\ref{thmproof} that $\Fs$ is contained in the subspace
$\{ y \in \Rs^{\mm}: y_{kl}=0\}$. Therefore, it follows from Theorem \ref{thmul} that
the pair of non-adjacent vertices $k$ and $l$ is universally linked.
\epr

\noindent{\bf Proof of Theorem \ref{thmmain2}}
Let $\Fs$ be the Cayley configuration space of framework $(G,p)$ and
suppose that Condition \ref{eqnmain2} of Theorem \ref{thmmain2} holds. Then it follows from
Theorem~\ref{thmproof} that $\Fs$ is contained in the trivial subspace
$\{ y \in \Rs^{\mm}: y=\bz \}$. Hence, $\Fs=\{\bz\}$ and hence, it follows from
from Theorem \ref{thmurF0} that $(G,p)$ is universally rigid.
\epr

Finally, the following corollary shows that Theorem \ref{thmsuff1} is a special case
of Theorem~\ref{thmmain2} when rank $\Omega = n-r-1$.

\begin{cor} \label{corend}
Let $(G,p)$ be an $r$-dimensional framework on $n$ vertices and let $\Omega$ be a positive semidefinite stress matrix
of $(G,p)$. If rank $\Omega = n-r-1$, then Condition (\ref{eqnmain2}) of Theorem \ref{thmmain2} is equivalent to
the statement that $(G,p)$ has no affine motions.
\end{cor}
\bpr Let $Z$ be a Gale matrix of $(G,p)$. Then the column space of $\Omega$ is equal to the column space of $Z$.
Hence, if Condition (\ref{eqnmain22}) holds, then the second statement of Lemma \ref{lemaffm} holds.
On the other hand, if  the second statement of Lemma \ref{lemaffm} holds, then
by Lemma \ref{lemxi}, Condition (\ref{eqnmain22}) holds.
\epr


\bibliographystyle{plain}
\end{document}